\documentclass{article}[12pt]
\usepackage{graphicx,pstricks,comment}
\usepackage[all]{xy}
\usepackage{subfigure}
\usepackage{rotating}
\usepackage{amsfonts}
\usepackage{amssymb}
\usepackage{amsmath}

\usepackage[utf8]{inputenc}

\usepackage{hyperref}

\usepackage{indentfirst}
\usepackage{epsfig}
\usepackage{psfrag}
\usepackage{enumerate}
\usepackage{array}
\usepackage{theorem}

\newcommand{\C}{\mathbb{C}}
\newcommand{\R}{\mathbb{R}}

\newcommand{\Z}{\mathbb{Z}}

\newcommand{\Q}{\mathbb{Q}}

\newcommand{\quat}{\mathbb{H}}

\newcommand{\la}{\langle}
\newcommand{\ra}{\rangle}

\newcommand{\nsub}{\mathrel{\unlhd}}

\newtheorem{thm}{Theorem}%[section]

\newtheorem{lem}{Lemma}%[section]
\newtheorem{prop}{Proposition}%[section]
{\theorembodyfont{\rmfamily} }
{\theorembodyfont{\rmfamily} }

\newcommand{\Pf}{{\em Proof}. }
\newcommand{\EPf}{\hfill$\Box$\vspace{.5cm}}
\title{Picard modular groups generated by complex reflections}
\author{Alice Mark, Julien Paupert\footnote{Second author partially supported by National Science Foundation Grant DMS-1708463.}, David Polletta}
\begin{document}
\maketitle

\begin{abstract} In this short note we use the presentations found in \cite{MP} and \cite{Po} to show that the 
Picard modular groups ${\rm PU}(2,1,\mathcal{O}_d)$ with $d=1,3,7$ (respectively the quaternion hyperbolic lattice ${\rm PSp}(2,1,\mathcal{H})$ with entries in the Hurwitz integer ring $\mathcal{H}$) are generated by complex (resp. quaternionic) reflections, and that the Picard modular groups ${\rm PU}(2,1,\mathcal{O}_d)$ with $d=2,11$ have an index 4 subgroup generated by complex reflections.
\end{abstract}

%\tableofcontents

\section{Introduction}

Hyperbolic reflection groups are an important class of groups in the realm of discrete subgroups and lattices in Lie groups, and more generally of discrete groups in geometry and topology. Such groups are accessible to a direct geometric description and understanding which are not always clear for groups defined algebraically or arithmetically.
While these reflection groups are relatively well understood in the constant curvature setting (they are then \emph{Coxeter groups} in Euclidean, spherical or real hyperbolic $n$-space), very little is known about their complex and quaternion hyperbolic counterparts.

 In the constant curvature setting, reflections are involutions whose fixed-point set is a totally geodesic submanifold of codimension 1. Such submanifolds do not exist in complex or quaternion hyperbolic space of dimension at least 2 (see \cite{CG}). Potential substitutes among isometries of complex hyperbolic $n$-space ${\rm H}_\C^n$ are \emph{complex
  reflections}, which are holomorphic isometries
fixing pointwise a totally geodesic complex hypersurface (a copy of
${\rm H}_\C^{n-1} \subset {\rm H}_\C^n$), and \emph{real reflections}, which are antiholomorphic involutions fixing
pointwise a Lagrangian subspace (a copy of ${\rm H}_\R^n \subset
{\rm H}_\C^n$). The situation in quaternion hyperbolic space is similar, with all totally geodesic subspaces being copies of  lower-dimensional real, complex or quaternion hyperbolic spaces.

Major open questions about the existence of lattices generated by real or complex reflections in ${\rm PU}(n,1)\simeq$ $ {\rm Isom}^0({\rm H}_\C^n)$ include the following. 
Do there exist lattices generated by (real or complex) reflections in ${\rm PU}(n,1)$ for all $n \geqslant 2$? For fixed $n \geqslant 2$, are there infinitely many (non-commensurable) lattices generated by (real or complex) reflections? In real hyperbolic space the answer to the first question is no by classical results of Vinberg (though sharp bounds on possible dimensions are far from known), whereas 2 and 3 are the only dimensions where infinitely many non-commensurable lattices generated by reflections are known to exist, by classical results of Poincar\'e and Andreev respectively.

In the complex hyperbolic case, Deligne--Mostow (\cite{DM}) and Mostow (\cite{Mos}) produced lattices generated by complex reflections in ${\rm PU}(n,1)$ for all $n \leqslant 9$ (finitely many for each $n$). Allcock found a further example in ${\rm PU}(13,1)$ in \cite{Al2}, related to the Leech lattice. In the quaternion case Allcock also produced in \cite{Al1} and \cite{Al2} lattices generated by quaternionic reflections in ${\rm PSp}(n,1)\simeq {\rm Isom}^0({\rm H}_\quat^n)$ in dimensions $n=2,3,5,7$, including the Hurwitz lattice in ${\rm PSp}(2,1)$ which we consider here (at least up to commensurability).

Stover showed that, among arithmetic lattices in ${\rm PU}(n,1)$ ($n \geqslant 2$), only those of \emph{first type} can contain complex reflections (Theorem 1.4 of \cite{St}, see also Example 9.2 of \cite{BFMS}). This class includes the Picard modular groups studied in this note. In particular there exist lattices in ${\rm PU}(n,1)$ for all $n \geqslant 2$ which do not contain a single complex reflection (even up to commensurability) - the arithmetic lattices of \emph{second type}, which in dimension 2 contains the fundamental groups of the so-called fake projective planes studied by Klingler (\cite{K}) and Prasad--Yeung (\cite{PY}) and classified by Cartwright--Steger (\cite{CS}). At the other extreme, it turns out that all known non-arithmetic lattices in ${\rm PU}(n,1)$ with $n \geqslant 2$ are commensurable to a lattice generated by complex reflections (there are 22 commensurability classes known for $n=2$ and two when $n=3$, see \cite{DPP} and \cite{D}).
%Arithmeticity, complex reflections, immersed totally geodesic surfaces....

The results in this note contribute to the small list of lattices in ${\rm PU}(2,1)$ known to be generated by complex reflections, among Picard modular groups with small discriminant. The Picard modular groups ${\rm PU}(2,1,\mathcal{O}_d)$ are the simplest kind of arithmetic lattices in ${\rm PU}(2,1)$, analogous to the Bianchi groups ${\rm PSL}(2,\mathcal{O}_d)$ in ${\rm PSL}(2,\C)$. (We denote by $\mathcal{O}_d$ the ring of integers of $\Q[\sqrt{-d}]$, where $d$ is a square-free positive integer). 
Bianchi proved in the seminal paper \cite{Bi} that the Bianchi groups are \emph{reflective}, i.e. generated by reflections up to finite index, for $d \leqslant 19$, $d \neq 14,17$. At the end of the 1980's, Shaiheev extended these results in \cite{Sh}, using results of Vinberg, proving that only finitely many of the Bianchi groups are reflective, including those with $d \leqslant 21$, $d \neq 14,17$. (The finiteness result now follows from a result of Agol, \cite{Ag}). The full classification of reflective Bianchi groups was obtained more recently in \cite{BeMc}.

The second author and Will proved in \cite{PWi} that the Picard modular groups ${\rm PU}(2,1,\mathcal{O}_d)$ are generated by real reflections when $d=1,2,3,7,11$. Our main result is the following:

\begin{thm}\label{main} The Picard modular groups ${\rm PU}(2,1,\mathcal{O}_d)$ with $d=1,3,7$ are generated by complex reflections; when $d=2,11$ they have an index 4 subgroup generated by complex reflections. The Hurwitz modular group ${\rm PSp}(2,1,\mathcal{H})$ is generated by quaternionic reflections.
\end{thm}

This result was known, at least up to finite index, for the Picard groups with $d=1,3$ and the Hurwitz group by work of Allcock \cite{Al1} (and later by \cite{FP} for $d=3$ and \cite{FFP} for $d=1$).

The authors would like to thank Matthew Stover for suggesting the direct computational method used in this note to determine indices of subgroups generated by reflections.

\section{Complex hyperbolic space and isometries}

We give a brief summary of key definitions and facts about complex hyperbolic space; see \cite{G} and \cite{CG} for more details (as well as \cite{KP} for quaternionic hyperbolic space, which we will not discuss in detail here as the aspects that we consider are similar to the complex case). We will consider only the case of dimension $n=2$ in this note, but the general setup is identical for higher dimensions so we state it for all $n \geqslant 1$. 
Consider $\C^{n,1}$, the vector space $\C^{n+1}$ endowed with a Hermitian form $\langle \cdot \, , \cdot \rangle$ of signature $(n,1)$.
Let $V^-=\left\lbrace Z \in \C^{n,1} | \langle Z , Z \rangle <0 \right\rbrace$.
Let $\pi: \C^{n+1}-\{0\} \longrightarrow \C{\rm P}^n$ denote projectivization.
Define ${\rm H}_\C^n$ to be $\pi(V^-) \subset \C{\rm P}^n$, endowed with the distance $d$ (Bergman metric) given by:

\begin{equation}\label{dist}
\cosh ^2 \frac{1}{2}d(\pi(X),\pi(Y)) = \frac{|\langle X, Y \rangle|^2}{\langle X, X \rangle  \langle Y, Y \rangle}
\end{equation}
From this formula it is clear that ${\rm PU}(n,1)$ acts by isometries on ${\rm H}_\C^n$ (where ${\rm U}(n,1)$ is the subgroup of ${\rm GL}(n+1,\C)$ preserving $\langle \cdot , \cdot \rangle$, and ${\rm PU}(n,1)$ is its image in ${\rm PGL}(n+1,\C)$). The boundary at infinity $\partial_\infty{\rm H}_\C^n$ is naturally identified with $\pi(V^0) \subset \C{\rm P}^n$, where $V^0=\left\lbrace Z \in \C^{n,1} | \langle Z , Z \rangle =0 \right\rbrace$.

\vspace{.3cm}

{\bf Fact:} ${\rm Isom}^0({\rm H}_\C^n)={\rm PU}(n,1)$, and ${\rm Isom}({\rm H}_\C^n)={\rm PU}(n,1) \ltimes \Z/2$ (complex conjugation).

\vspace{.3cm}

{\bf Classification:} $g \in {\rm PU}(n,1)$ is of one of the following types:
\begin{itemize}
\item \emph{elliptic}: $g$ has a fixed point in ${\rm H}_\C^n$
\item \emph{parabolic}: $g$ has (no fixed point in ${\rm H}_\C^n$ and) exactly one fixed point in $\partial_\infty{\rm H}_\C^n$
\item \emph{loxodromic}: $g$ has (no fixed point in ${\rm H}_\C^n$ and) exactly two fixed points in $\partial_\infty{\rm H}_\C^n$
 \end{itemize}

{\bf Definitions:} For any $1\leqslant k \leqslant n$, a \emph{complex k-plane} is a $k$-dimensional projective subspace of $\C P^n$ intersecting $\pi(V^-)$ non-trivially (so, it is an isometrically embedded copy of ${\rm H}_\C^{k} \subset {\rm H}_\C^n$). Complex 1-planes are usually called \emph{complex lines}. A \emph{complex reflection} is an elliptic isometry $g\in {\rm PU}(n,1)$ whose fixed-point set is a complex $(n-1)$-plane. 
%An elliptic isometry $g$ is called \emph{regular} if any of its matrix representatives $A \in U(n,1)$ has distinct eigenvalues. 
The eigenvalues of a matrix $A \in {\rm U}(n,1)$ representing an elliptic isometry $g$ have modulus one. Exactly one of these eigenvalues has eigenvectors in $V^-$ (projecting to a fixed point of $g$ in ${\rm H}_\C^n$), and such an eigenvalue will be called \emph{of negative type}. An elliptic isometry $g\in {\rm PU}(n,1)$ is a complex reflection if and only if the negative type eigenvalue of any of its matrix representatives has multiplicity $n$.
%Regular elliptic isometries have an isolated fixed point in ${\rm H}_\C^n$.   

\section{Picard modular groups and complex reflections}

We use the \emph{Siegel model} of hyperbolic space ${\rm H}_\C^2$, which is the projective model associated to the Hermitian form on $\C^{3}$ given by $\langle Z,W \rangle = W^* J Z$ with:

$$ J=\left(\begin{array}{ccc}
0 & 0 & 1 \\
0 & 1 & 0 \\
1 & 0 & 0 \end{array}\right)
$$

The \emph{Picard modular groups} are the arithmetic lattices $\Gamma_d={\rm PU}(2,1,\mathcal{O}_d)$ in ${\rm PU}(2,1)$, where $d$ is a squarefree positive integer and $\mathcal{O}_d$ is the ring of integers of $\Q[i\sqrt{d}]$. The following elements belong to $\Gamma_d$ for all $d$ (taking $u=i$ when $d=1$, $u=e^{2i\pi/3}$ when $d=3$, and $u=-1$ for all other values of $d$): 
$$
\begin{array}{cc} I_0=\left[ \begin{array}{ccc} 0 & 0 & 1 \\ 0 & -1 & 0 \\ 1 & 0 & 0 \end{array}\right],  \ & \ 
  R=\left[ \begin{array}{ccc} 1 & 0 & 0 \\ 0 & u & 0 \\ 0 & 0 & 1 \end{array}\right].  
  \end{array}
$$

\begin{lem} $I_0$ and $R$ are complex reflections.  
\end{lem}

\Pf $R$ visibly has eigenvalues $\{1,1,u \}$ and 1 is of negative type, hence $R$ is a complex reflection. Likewise, $I_0$ has eigenvalues $\{-1,-1,1 \}$ with -1 of negative type (its eigenspace is the span of $e_2$ and $e_1-e_3$), hence $I_0$ is a complex reflection. Note that $I_0$ has order 2, and $R$ has order 2 except when $d=3$ (when it has order 6) and $d=1$ (when it has order 4). \EPf 

Given elements $\gamma_1,...,\gamma_k$ of a group $\Gamma$, we denote $\la \la \gamma_1,...,\gamma_k \ra \ra$ the normal closure of $\gamma_1,...,\gamma_k$ in $\Gamma$, that is the smallest normal subgroup of $\Gamma$ containing $\gamma_1,...,\gamma_k$.
 
\begin{prop} Denote as above $\Gamma_d={\rm PU}(2,1,\mathcal{O}_d)$.
\begin{itemize}
\item When $d=3$, $\la \la R \ra \ra = \Gamma_3$.
\item When $d=1$, $\la \la R, I_0 \ra \ra = \Gamma_1$,  $\la \la I_0 \ra \ra$ has index 4 in $\Gamma_1$ and $\la \la R  \ra \ra$ has index 96 in $\Gamma_1$.
\item When $d=7$, $\la \la I_0 \ra \ra = \Gamma_7$ and $\la \la R  \ra \ra$ has index 168 in $\Gamma_7$.
\item When $d=2$, $\la \la R, I_0 \ra \ra = \la \la I_0 \ra \ra $ has index 4 in $\Gamma_2$.
\item When $d=11$, $\la \la R, I_0 \ra \ra$ has index 4 in $\Gamma_{11}$ and $\la \la R \ra \ra$ has index 13,824 in $\Gamma_{11}$.

\end{itemize}
\end{prop}

\Pf Let $\Gamma$ be a group with (say, finite) presentation $\la s_1,...,s_n  \, | \, r_1,...,r_p \ra$ and $w_1,...,w_k$ elements of $\Gamma$ given as words in the generators $s_1,...,s_n$. Then, by construction of group presentations:

$$ \hat{\Gamma}= \la s_1,...,s_n  \, | \, r_1,...,r_p, w_1,...,w_k \ra \simeq \Gamma / \la \la w_1,...,w_k \ra \ra.
$$ 
In particular, the order of $ \hat{\Gamma}$ is equal to the index of $\la \la w_1,...,w_k \ra \ra$ in $\Gamma$. The statements in the proposition give the results of this procedure applied to the presentations for the $\Gamma_d$ obtained in \cite{MP} and \cite{Po}, available as Magma files at \cite{MCode} and \cite{PoCode}. (Note that both $R$ and $I_0$ appear conveniently as generators in these presentations.) More specifically, we add the relation $R$ (resp. $I_0$, resp. $R$ and $I_0$) to these presentations and use the Magma command ${\tt Order(G);}$ to compute the order of the quotient. Most of these computations can be done in a matter of seconds with the online Magma Calculator available at \cite{Mag} (except for the higher orders 96, 168 and 13,824 which require an installed version of Magma). \EPf

\section{The Hurwitz modular group and quaternionic reflections}

Quaternionic hyperbolic space admits a Siegel model analogous to the one discussed above for complex hyperbolic space, with the usual caveats of linear algebra over the quaternions (eg. matrices act on vectors by left multiplication and scalars by right multiplication). See \cite{KP} or \cite{MP} for details. A \emph{quaternionic reflection} is an elliptic isometry $g\in {\rm PSp}(n,1)$ whose fixed-point set is a quaternionic $(n-1)$-plane. We are interested in the lattice $\Gamma_\mathcal{H}={\rm PSp}(2,1,\mathcal{H})<{\rm PSp}(n,1)$, consisting of the (projectivized) matrices in ${\rm PSp}(n,1)$ whose entries lie in the Hurwitz ring $\mathcal{H}=\Z[1,i,j,k,\sigma] \subset \quat$, denoting $\sigma=\frac{1+i+j+k}{2}$. The relevant elements of $\Gamma_\mathcal{H}$ for our purposes are:
$$
\begin{array}{ccc} I_0=\left[ \begin{array}{ccc} 0 & 0 & 1 \\ 0 & -1 & 0 \\ 1 & 0 & 0 \end{array}\right],  \ & \ 
  R_i=\left[ \begin{array}{ccc} 1 & 0 & 0 \\ 0 & i & 0 \\ 0 & 0 & 1 \end{array}\right],  \ & \ 
  R_\sigma=\left[ \begin{array}{ccc} 1 & 0 & 0 \\ 0 & \sigma & 0 \\ 0 & 0 & 1 \end{array}\right]. 
\end{array}
$$

The same computation as the proof of Lemma 1 gives the following:   
 \begin{lem} $I_0$, $R_i$ and $R_\sigma$ are quaternionic reflections.  
\end{lem}
  
\begin{prop} Let $\Gamma_\mathcal{H}$ and $I_0, R_i, R_\sigma \in \Gamma_\mathcal{H}$ as above. Then $\la \la R_\sigma \ra \ra = \Gamma_\mathcal{H}$, $\la \la I_0  \ra \ra$ has index dividing 12 in $\Gamma_\mathcal{H}$ and 
$\la \la R_i  \ra \ra$ has index dividing 648 in $\Gamma_\mathcal{H}$. 
\end{prop}

\Pf The procedure is the same as for the proof of Proposition 1, the difference being that we only use a partial presentation for $\Gamma_\mathcal{H}$, that is a presentation with all generators but only some of the relations. As observed in \cite{MP}, the presentation obtained there for $\Gamma_\mathcal{H}$ is too large for Magma to handle directly (it has 33 generators and 968,480 relations, and the text file is a bit over 200 MB). Rather, we use the partial presentation obtained by keeping only the 1000 first relations -- a text file for this presentation is available as {\tt QuaternionsTruncated1000.txt} at \cite{MCode}. 

Denoting $\tilde{\Gamma}$ the abstract group with this partial presentation, we have a surjective homomorphism $\pi:\tilde{\Gamma} \longrightarrow \Gamma_\mathcal{H}$, obtained by adding the remaining relations for $\Gamma_\mathcal{H}$. Given a normal subgroup $H \nsub \tilde{\Gamma}$, $\pi$ induces a surjective homomorphism $\tilde{\Gamma}/H \longrightarrow \Gamma_\mathcal{H}/\pi(H)$. Taking $H=\la \la R_\sigma \ra \ra, \, \la \la I_0  \ra \ra, \, \la \la R_i  \ra \ra$ successively, we compute the order of $\tilde{\Gamma}/H$ as above by adding the single relation $R_\sigma, I_0, R_i$ to the presentation for $G=\tilde{\Gamma}$ using the Magma command {\tt Order(G);} and the result follows. \EPf

\raggedright

\vspace{.2cm}

\begin{flushleft}
 \textsc{Alice Mark\\
   Department of Mathematics, Vanderbilt University}\\
       \verb|alice.h.mark@vanderbilt.edu| 
\end{flushleft}

\vspace{.2cm}

\begin{flushleft}
  \textsc{Julien Paupert, David Polletta\\
   School of Mathematical and Statistical Sciences, Arizona State University}\\
       \verb|paupert@asu.edu|, \verb|dpollett@asu.edu|
\end{flushleft}

\end{document}